\documentclass{amsart}
\usepackage{amsmath}

\title{A novel count of the spanning trees of a cube}

\author{Thomas W.\ Mattman}
\address{Department of Mathematics and Statistics,
California State University, Chico,
Chico, CA 95929-0525}
\email{TMattman@CSUChico.edu}


\begin{document}

\begin{abstract}
Using the special value at $u=1$ of the Artin-Ihara $L$-function, we give a short 
proof of the count of the number of spanning trees in the $n$-cube. 
\end{abstract}

\maketitle

The special value at $u=1$ of the Artin-Ihara $L$-function is
an innovative way to count spanning trees in graph covers.
For example, in~\cite{HMSV}, we made two observations about
$\kappa_G$, the number of spanning trees in graph $G$.
First, if $G$ is an abelian cover of $H$, then $\kappa_H \mid \kappa_G$.
Note that Baker and Norine~\cite{BN} also prove this and, not just for
abelian covers, but also for every topological cover, and even more generally, 
for every non-constant harmonic morphism of graphs.
    
For the second observation, let $G$ be a $(\mathbb{Z}/2 \mathbb{Z})^m$ cover of $H$ 
and $H_1, H_2, H_3, \ldots , H_{2^m-1}$, the intermediate double covers. Then,
\begin{equation}
\kappa_{G} = \frac{2^{2^{m}-m-1}}{\kappa_{H}^{2^{m} - 2}}\prod_{i=1}^{2^{m}-1}\kappa_{H_i}. 
\label{eqZm}
\end{equation}
Moreover, we speculate that this type of relationship generalizes to other abelian covers. 

To encourage study of such generalizations, we illustrate how this formula leads
to a novel count of the number of trees in the $n$-cube, $C_n$.
We were inspired by Stanley's book \cite{S}, which uses eigenvalues to deduce
\begin{equation}
\kappa_{C_n} = 2^{2^{n}-n-1}\prod_{i=1}^{n} i^{\binom{n}{i}}.
\label{eqkcn}
\end{equation}
Aside from the striking similarity of the two equations,
we were intrigued that an early version of the text remarked on the lack of a combinatorial
proof of Equation~\ref{eqkcn}. In the meantime, Bernardi~\cite{B} and Biane and Chapuy~\cite{BC}
have provided combinatorial arguments. Our new 
proof is likewise combinatorial in flavour and also similar to
the other two in that we do not have an explicit bijection.

For $n \geq 1$, let $C_n$ denote the $n$-cube: the vertices are the $2^n$ points in $(\mathbb{Z}/2 \mathbb{Z})^n$ 
and two vertices are adjacent if they differ at exactly one coordinate.  Our base graph $H$,
a generalized theta graph, is the multigraph
with two vertices $x$ and $y$ and $n$ edges between them. This means the Galois group of $C_n$ over $H$ is
$(\mathbb{Z}/2 \mathbb{Z})^{n-1}$. We can describe the double covers as examples of $B_{a,n-a}$ graphs.
A $B_{a,b}$ is bipartite with four vertices, $\{x_1,x_2,y_1,y_2\}$, the $x_i$'s forming one part and the 
$y_i$'s the other. There are $a$ $x_iy_i$ edges for $i = 1,2$ (making $2a$ in total) and $b$ edges
of the form $x_1y_2$ and a further $b$ $x_2y_1$ edges. 
Thus, $B_{a,b}$ is $(a+b)$-regular and has $2(a+b)$ edges.
Note that for $1 \leq a \leq n-1$, 
$B_{a,n-a}$ is a double cover of $H$ with Galois group $\mathbb{Z}/2 \mathbb{Z}$.

For convenience, we assume that $n \geq 3$ is odd. The argument for even $n$ is similar.
Observe that there are $2^{n-1}-1$ double covers $H_i$ of the form $B_{a,n-a}$. Indeed,
for each $1 \leq a \leq (n-1)/2$, there are $\binom{n}{a}$ $B_{a,n-a}$ graphs that correspond
to choosing $a$ of the $n$ edges of $H$. This gives a total of
$$\sum_{a = 1}^{(n-1)/2} \binom{n}{a} = 2^{n-1}-1$$
graphs. The final ingredients in our computation are the observations 
that $\kappa_{H} = n$ and $\kappa_{B_{a,n-a}} = 2a(n-a)n$.

Applying Equation~\ref{eqZm} with $m = n-1$ we have
\begin{eqnarray*}
K_{C_n} & = & \frac{2^{2^{n-1}-(n-1)-1}}{\kappa_{H}^{2^{n-1} - 2}}\prod_{i=1}^{2^{n-1}-1}\kappa_{H_i} \\
& = & \frac{2^{2^{n-1}-n}}{n^{2^{n-1} - 2}}\prod_{a=1}^{(n-1)/2}(\kappa_{B_{a,n-a}})^{\binom{n}{a}} \\
& = & \frac{2^{2^{n-1}-n}}{n^{2^{n-1} - 2}}\prod_{a=1}^{(n-1)/2}(2a(n-a)n)^{\binom{n}{a}} \\
& = & \frac{2^{2^{n-1}-n}}{n^{2^{n-1} - 2}}(2n)^{\sum_{a=1}^{(n-1)/2} \binom{n}{a}}\prod_{a=1}^{(n-1)/2}a^{\binom{n}{a}}(n-a)^{\binom{n}{n-a}} \\
& = & \frac{2^{2^{n-1}-n}}{n^{2^{n-1} - 2}}(2n)^{2^{n-1}-1} \prod_{a=1}^{(n-1)}a^{\binom{n}{a}} \\
& = & 2^{2^{n}-n-1}\prod_{a=1}^{n} a^{\binom{n}{a}},
\end{eqnarray*}
recovering Equation \ref{eqkcn}.

\medskip

\noindent%
{\bf Acknowledgement:} We thank Kyle Hammer, Jonathan Sands, Daniel Valli\`eres, and Matt Baker for helpful conversations.

\end{document}